\documentclass[12pt]{amsart}
\usepackage{hyperref}
\usepackage{amsfonts}
\usepackage{amssymb}
\def\depth{\operatorname{depth}}
\def\gr{\operatorname{gr}}
\newcommand{\M}{\mathcal M}
\newcommand{\m}{\mathfrak m}

\newtheorem{theorem}{Theorem}
\newtheorem{lemma}[theorem]{Lemma}

\newtheorem{proposition}[theorem]{Proposition}
\newtheorem{remark}[theorem]{Remark}
\newtheorem{example}[theorem]{Example}

\begin{document}
\title[On the depth of graded rings associated to lex-segment ideals]
{On the depth of graded rings associated to lex-segment
ideals in $K[x,y]$}
\author{A. V. Jayanthan}
\address{Department of Mathematics, Indian Institute of Technology
Madras, Chennai 600036, INDIA.}

\begin{abstract}
In this article, we show that the depths of the associated graded ring
and fiber cone of a lex-segment ideal in $K[x,y]$ are equal.
\end{abstract}

\maketitle
\vskip 4mm
\noindent

\section{Introduction}

Let $K$ be a field of characteristic zero and $R = K[x_1, \ldots,
x_n]$ be the polynomial ring in $n$ variables over $K$. Let $R_i$
denote the $K$-vector subspace of all monomials of degree $i$. We fix
the ordering of variables as $x_1 > x_2 > \cdots > x_n$. For
monomials $u = x_1^{a_1}\cdots x_n^{a_n}$ and $v = x_1^{b_1}\cdots
x_n^{b_n}$, we say that $u <_{Lex} v$ if $\deg u \leq \deg v$ or $\deg
u = \deg v$ and $b_i - a_i > 0$ for the first time when it is nonzero.
An initial lex-segment in degree $d$ is the set of all monomials of
the form $\{m \in R_d \; : \; m \geq u \}$, where $u \in R_d$. A
graded ideal $I$ is said to be a {\it lex-segment ideal} if $I_d$ is
generated by initial lex-segments for each $d$ with $I_d \neq 0$.
Lex-segment ideals are important due to many reasons.
It is well known that among ideals with a given Hilbert function, the
lex-segment ideal has the largest number of generators. A. M. Bigatti
\cite{big}
and H. A. Hulett \cite{hul} in characteristic zero and K. Pardue
\cite{p} in positive characteristic
generalized this to all Betti numbers. They proved that the
lex-segment ideals have the largest Betti numbers among all ideals
with a given Hilbert function. Lex-segment ideals are of interest also
due to classical reasons. O. Zariski used the theory of contracted
ideals to study complete ideals in $2$-dimensional regular local rings
$(R,\m)$. In the graded setting, when $K$ is algebraically closed,
Zariski's factorization theorem for homogeneous contracted ideals
asserts that any homogeneous contracted ideal $I$ can be written as $I
= \m^c L_1\cdots L_t$, where each $L_i$ is a lex-segment ideal with
respect to an appropriate system of coordinates $x_i, y_i$ which
depends on $i$. \cite[Theorem 1, Appendix 5]{zs}, \cite[Theorem
3.8]{cdjr}.

In this article, we study the blowup algebras, namely, the associated
graded ring and the fiber cone of lex-segment ideals in a two dimensional
polynomial ring. Let $R$ be a ring, $I$ any ideal of $R$ and $\m$ a
maximal ideal. Then the associated graded ring and the fiber cone of
$I$ are respectively defined as $\gr_I(R) = \oplus_{n\geq
0}I^n/I^{n+1}$ and $F(I) = \oplus_{n\geq 0}I^n/\m I^n$.
In \cite{hm}, Huckaba and Marley showed that in a regular local ring
$(R,\m)$, $\depth \gr_I(R) = \depth R(I) - 1$ for any $\m$-primary
ideal $I$, where $R(I) = \oplus_{n\geq0}I^nt^n$ denotes the Rees
algebra of $I$. It is interesting to ask if there is a similar relation
between the depths of the fiber cone and the associated graded ring.
It is well known that this is not the case in general (cf. Example
\ref{ex1}, Example \ref{ex2}). In this article, we prove that the
depths of these algebras are equal for lex-segment ideals in $K[x,y]$,
where $K$ is a field of characteristic zero.

\vskip 2mm
\noindent
\textbf{Acknowledgements:} The author would like to thank Aldo Conca,
M. E. Rossi, G. Valla, J. K. Verma and S. Goto for useful discussions
regarding the contents of the paper.

\section{Equality of depths}
Let $R = K[x,y]$, where $K$ is a field of characteristic zero and $\M
= (x,y)$. In this case, the lex-segment ideals are easy to describe.
If $I$ is a lex-segment ideal in $K[x,y]$, then $I = (x^d,
x^{d-1}y^{a_1}, \ldots, x^{d-k}y^{a_k})$ for some $0 \leq k \leq d$
and $1\leq a_1 < a_2 < \cdots < a_k$. Note that if $I$ is a
lex-segment ideal, then $I^n$ is also a lex-segment ideal for all $n
\geq 1$.

\begin{remark}\label{localgl}
Let $S = K[\![x,y]\!]$ and $\m = (x,y)$. Then for any ideal $I \subset
\m$, $S/IS \cong R/I$ and $S/\m S \cong R/\m$. Therefore $\gr_{IS}(S)
\cong \gr_I(R)$ and $F(IS) \cong F(I)$ \cite[Lemma 2.1]{cdjr}.
Hence, we may use the local techniques to prove the results for
$\gr_IS(S)$ and $F(IS)$ and derive the same for $\gr_I(R)$ and $F(I)$.
\end{remark}

We first show that the Cohen-Macaulay property of the associated
graded ring and the fiber cone are equivalent. The dimension of the
fiber cone, denoted by $s(I)$, is called the analytic spread. It is
well known that $h(I) \leq s(I)$, where $h(I)$ denote the height of
the ideal $I$. The difference, $s(I) - h(I),$ is called the analytic
deviation. Let $I = (x^d, x^{d-1}y^{a_1}, \ldots, x^{d-k}y^{a_k})$. If
$k = d$, then $I$ is an $\M$-primary homogeneous contracted ideal.
Because of Remark \ref{localgl}, we can use local theory of
$\M$-primary contracted ideals in 2-dimensional regular local rings to
study the blowup algebras. If $0 < k < d$, then $I$ is a
non-$\M$-primary ideal of analytic deviation one. 
Here we note that if $0 < k < d$, then $I = x^{d-k}(x^k,
x^{k-1}y^{a_1}, \ldots, xy^{a_{k-1}}, y^{a_k})$, which is of the
form $I = zL$, where $z$ is an $R$-regular element and $L$ an
$\M$-primary homogeneous contracted ideal. 
We show that the depth of $\gr_I(R)$ is at most
the depth of $\gr_L(R)$. In particular, when $\gr_I(R)$ is
Cohen-Macaulay, so is $\gr_L(R)$. For an element $a \in I$, let $a^*$
denote its initial form in $\gr_I(R)$ and $a^o$ denote its initial
form in $F(I)$.

Let $I$ be an ideal of a ring $R$. An ideal $J \subseteq I$ is said to
be reduction of $I$ if $I^{n+1} = JI^n$ for some $n \geq 0$. A
reduction which is minimal with respect to inclusion is called a
minimal reduction. For a reduction $J$ of $I$, the number $r_J(I) =
\min\{n \; |\; I^{n+1} = JI^n\},$ is called the reduction number of
$I$ with respect to $J$.

\begin{proposition}\label{grdepth}
Let $(R,\m)$ be a Noetherian local ring and $L$ an $\m$-primary ideal
of $R$. Let $x$ be a regular element in $R$ and $I = xL$. Then $\depth
\gr_I(R) \leq \depth \gr_L(R)$. In particular, if $\gr_I(R)$ is
Cohen-Macaulay, then so is $\gr_L(R)$.
\end{proposition}
\begin{proof}
Let $\depth \gr_I(R) = t$. Let $a_1, \ldots, a_t \in L \; \backslash
\; L^2$ and $b_i = xa_i$ be such that $b_1^*, \ldots, b_t^* \in
\gr_I(R)$ is a regular sequence. Then by Valabrega-Valla \cite{vv},
$(b_1, \ldots, b_t) \cap I^n = (b_1, \ldots, b_t)I^{n-1}$ for all $n
\geq 1$. We show that $(a_1, \ldots, a_t) \cap L^n =
(a_1, \ldots, a_t)L^{n-1}$ for all $n \geq 1$.

Let $p \in (a_1, \ldots, a_t) \cap L^n$ for some $n \geq 1$. Then
$x^np \in (b_1, \ldots, b_t) \cap I^n = (b_1, \ldots, b_t)I^{n-1} =
x^n (a_1, \ldots, a_t)L^{n-1}$. Therefore $x^np = x^nq$ for some $q
\in (a_1, \ldots, a_t)L^{n-1}$. Since $x$ is regular in $R$, $p = q$
which implies that $p \in (a_1, \ldots, a_t)L^{n-1}$. Therefore, by
Valabrega-Valla condition, $a_1^*, \ldots, a_t^* \in \gr_L(R)$ is a
regular sequence.
\end{proof}

\begin{remark}
In the above Proposition, we have shown that if $b_1^*, \ldots, b_t^*$
is a regular sequence in $\gr_I(R)$, then $a_1^*, \ldots, a_t^*$ is a
regular sequence in $\gr_L(R)$. The following example shows that the
converse is not true in general.
\end{remark}

\begin{example}
Let $R = K[x,y]$. Let $L = \M = (x, y)$ and $I = (x^3, x^2y)$. Then
$x^*, y^*$ is a regular sequence in $\gr_L(R)$. It can be easily seen
that $I^2 : x^3 = (x^3, x^2y, xy^2) \neq I$. Therefore $(x^3)^* \in
\gr_I(R)$ is not regular. However, this does not imply that the
$\depth \gr_I(R) < 2$. In fact, in this case, it can be seen (using
any of the computational commutative algebra packages) that
$\gr_I(R)$ is indeed Cohen-Macaulay.
\end{example}

The following result follows directly from Theorem 2.1 of \cite{drv}.

\begin{proposition}\label{fibercm}
Let $(R,\m)$ be a Cohen-Macaulay local ring and $I$ be an ideal of $R$
with $s(I) = r$ and
$$
H(F(I),t) = \frac{a+bt}{(1-t)^r}.
$$
If $F(I)$ is Cohen-Macaulay, then $r_J(I) \leq 1$ for any minimal
reduction $J$ of $I$.
\end{proposition}

We show that the Cohen-Macaulay property of the associated graded ring
and the fiber cone are equivalent:

\begin{theorem}
Let $I$ be a lex-segment ideal in $K[x,y]$. Then $F(I)$ is
Cohen-Macaulay if and only if $\gr_I(R)$ is Cohen-Macaulay.
\end{theorem}

\begin{proof}
Let $I = (x^d, x^{d-1}y^{a_1}, \ldots, x^{d-k}y^{a_k})$ for some $0
\leq k \leq d$ and $1 \leq a_1 < a_2 < \cdots < a_k$. If $k = 0$, then
$I = (x^d)$ and both $\gr_I(R)$ and $F(I)$ are Cohen-Macaulay. We deal
the cases $k = d$ and $0 < k < d$ separately. Note that because
of Remark \ref{localgl}, we may assume that $I$ is an ideal in a two
dimensional regular local ring $(R,\M)$.

\vskip 2mm
\noindent
Let $k = d$. In this case, $I$ is $\M$-primary. Suppose $\gr_I(R)$ is
Cohen-Macaulay. Since $I$ is contracted, by Theorem 5.1 of \cite{h},
for any minimal reduction $J \subset I, \; I^2 = JI$. By \cite{s},
$F(I)$ is Cohen-Macaulay.
\vskip 2mm
\noindent
Conversely, suppose that $F(I)$ is Cohen-Macaulay. Note that for all
$n \geq 0$, $\mu(I^n) = nd + 1$ so that the Hilbert series of $F(I)$
is given by
$$
H(F(I), t) = \frac{1 + (d-1)t}{(1-t)^2}.
$$
Since $F(I)$ is Cohen-Macaulay, by Proposition \ref{fibercm}, $r_J(I)
\leq 1$ for any minimal reduction $J$ of $I$. Therefore, $\gr_I(R)$ is
Cohen-Macaulay \cite{vv}.
\vskip 2mm
\noindent
Now let $0 < k < d$. In this case, $I = x^{d-k}L$, where $L =
(x^k, x^{k-1}y^{a_1}, \ldots, y^{a_k})$.
Suppose $\gr_I(R)$ is Cohen-Macaulay. Then by Proposition
\ref{grdepth}, $\gr_L(R)$ is Cohen-Macaulay. By Proposition 2.6 of
\cite{hm}, $R(L)$ is Cohen-Macaulay. Hence the reduction number $r(L)$
is at most one, by Goto-Shimoda theorem \cite{gs}. Therefore $r(I) \leq 1$.
Therefore by \cite{s}, $F(I)$ is Cohen-Macaulay.

\vskip 2mm
\noindent
Suppose now that $F(I)$ is Cohen-Macaulay. Since $\mu(I^n) = nk + 1 $,
$$
H(F(I),t) = \frac{1+(k-1)t}{(1-t)^2}.
$$
Therefore by Proposition
\ref{fibercm}, $I^2 = JI$ for any minimal reduction $J$ of $I$.
Hence, Valabrega-Valla condition implies that $\depth \gr_I(R)
\geq s(I) = 2$ so that $\gr_I(R)$ is Cohen-Macaulay.
\end{proof}

Using Proposition \ref{grdepth}, we give a simple proof of the fact
that for lex-segment ideals the Cohen-Macaulayness of the Rees algebra
and the associated graded rings are equivalent. This has been proved
for $\m$-primary ideals in a regular local ring $(R,\m)$. Since we
could not find a generalization of this result for the
non-$\m$-primary ideals, we use this opportunity to present a simple
proof in the case of lex-segment ideals.

\begin{theorem}
Let $R = K[x,y]$ and $I$ a lex-segment ideal. Then $R(I)$ is
Cohen-Macaulay if and only if $\gr_I(R)$ is Cohen-Macaulay.
\end{theorem}

\begin{proof}
Let $I = (x^d, x^{d-1}y^{a_1}, \ldots, x^{d-k}y^{a_k})$. If $k = d$,
then $I$ is $\M$-primary and hence it follows from Proposition 2.6 of
\cite{hm}. If $k = 0$, then $I$ is a parameter ideal and hence both
the graded algebras are Cohen-Macaulay. Suppose $0 < k < d$. Then $I =
x^{d-k}L$, where $L = (x^k, x^{k-1}y^{a_1}, \ldots, y^{a_k})$. If
$\gr_I(R)$ is Cohen-Macaulay, then by Proposition \ref{grdepth}
$\gr_L(R)$ is Cohen-Macaulay. Since $L$ is $\M$-primary by Proposition
2.6 of \cite{hm}, $R(L)$ is Cohen-Macaulay. Since $x^{d-k}$ is a
regular element, $R(L) \cong R(I)$ and hence $R(I)$ is Cohen-Macaulay.

\vskip 2mm
\noindent
Conversely, suppose $R(I)$ is Cohen-Macaulay. Hence $R(L)$
Cohen-Macaulay. By Goto-Shimoda theorem, $r(L) \leq 1$. Therefore
$r(I) \leq 1$ and hence $\gr_I(R)$ is Cohen-Macaulay.
\end{proof}

\begin{remark}
The above result together with Theorem 3.4 of \cite{m} implies that
for any lex-segment ideal $I$ in $K[x,y]$, $\depth \gr_I(R) = \depth
R(I) - 1$.
\end{remark}

Now we proceed to prove that the fiber cone has positive depth if and
only if the associated graded ring has positive depth. We begin with
some properties of lex-segment ideals.
\begin{lemma}
Let $I = (x^d, x^{d-1}y^{a_1}, \ldots, x^{d-k}y^{a_k})$ be a lex-segment
ideal in $R = K[x,y]$ and $\M = (x,y)$. Then,
\begin{enumerate}
\item $\M I^n : y = I^n$ for all $n \geq 0$.
\item $\M I^{n+1} : I = \M(I^{n+1} : I)$ for all $n \geq 0$.
\end{enumerate}
\end{lemma}

\begin{proof}
(1) Note that $\M I^n = (x^{nd+1}) + yI^n$ for all $n \geq 0$. Since
$\M I^n : y$ is a monomial ideal, it is enough to show that the
monomials in $\M I^n : y$ are in $I^n$. For a polynomial $p \in
K[x,y]$, let $\deg_x p$ denote the degree of the polynomial with
respect to $x$, considering it as a polynomial in $x$ with
coefficients in $K[y]$ and $\deg_y p$ denote the degree of the
polynomial $p$ with respect to $y$, considering it as a polynomial in
$y$ with coefficients in $K[x]$. Let $p \in \M I^n : y$. If $\deg_x p
\geq nd$, then clearly $p \in I^n$. Therefore, we may assume that
$\deg_x p < nd$. Set $p = x^{nd-r}y^s$ for some $r, s \geq 1$. Since
$I^n$ is also a lex-segment ideal, for each $nd-nk \leq t \leq nd$,
there exists a unique minimal generator $p_t$ such that $\deg_x p_t =
t$.  Let $u = x^{nd-r}y^b$ be the minimal generator of $I^n$
with $\deg_x u = nd-r$. Then, $py \in \M I^n$ implies that $s+1 \geq
b + 1$. Hence $s \geq b$. Therefore, $p = x^{nd-r}y^s
\in I^n$.

\vskip 2mm
\noindent
(2) Let $p \in \M I^{n+1} : I$. If $p = x^r$ for some $r$, then
$x^r.x^d \in \M I^{n+1}$. Since any term which is a pure power in $x$
in $\M I^{n+1}$ has degree at least $(n+1)d+1$, we get that $r+d \geq
(n+1)d+1$. Therefore $r \geq nd+1$ so that $x^r \in \M I^n \subseteq
\M (I^{n+1} : I)$. Now assume that $y$ divides $p$. Write $p = yp'$.
Then $yp'f \in \M I^{n+1}$ for all $f \in I$. Therefore $p'f \in \M
I^{n+1} : y = I^{n+1}$ for all $f \in I$. Hence $p' \in I^{n+1} : I$
so that $p \in \M (I^{n+1} : I)$.
\end{proof}

\begin{theorem}
Let $I$ be a lex-segment ideal in $R$. Then $\depth \gr_I(R)
> 0$ if and only if $\depth F(I) > 0$.
\end{theorem}

\begin{proof}
Let $\depth \gr_I(R) > 0$. Then, $I^{n+1} : I = I^n$ for all $n \geq
0$. Therefore, $\M I^{n+1} : I = \M (I^{n+1} : I) = \M I^n$ for all $n
\geq 0$. Hence $\depth F(I) > 0$.
\vskip 2mm
\noindent
Conversely, assume that $\depth F(I) > 0$. Then, $(\M I^{n+1} : I)
\cap I^n = \M I^n$ for all $n \geq 0$. We need to prove that $I^{n+1}
: I = I^n$ for all $n \geq 0$. Suppose that there exists an $n$ such
that $I^n \varsubsetneq I^{n+1} : I$. Since $I$ is a monomial ideal,
$I^{n+1} : I$ is generated by monomials and hence there exists a
monomial generator $p$ of $I^{n+1} : I$ such that $p \notin I^n$.
Since $I^{n+1} : I$ is also a lex-segment ideal, we can write $p =
x^{nd-t}y^s$ for some $s$. Let $q \in I^n$ be the minimal generator of
$I^n$ such that $\deg_x q = nd-t$. Then $\deg_y q > s$, since $p
\notin I^n$. Therefore $qI \subseteq \M I^{n+1}$. Hence $q \in (\M
I^{n+1} : I) \cap I^n = \M I^n$. This contradicts the fact that $q$ is
a minimal generator of $I^n$. Therefore $I^{n+1} : I = I^n$ for all $n
\geq 0$. Hence $\depth \gr_I(R) > 0$.
\end{proof}

The result that $\depth \gr_I(R) = \depth F(I) = \depth R(I) - 1$ has
been proved for $\M$-primary lex-segment ideals in $k[x,y]$ by Conca,
De Negri and Rossi in \cite{cdr}.

\section{Examples}
In this section we give some examples to show that depths of fiber
cone and the associated graded rings are not related in general.
\begin{example}\label{ex1}
Let $I = (x^5, x^3y^3, xy^7, y^9) \subset R = K[x,y]$. Then it can be
seen that $x^2y^6 \in I^2 : I,$ but not in $I$. Therefore $\depth
\gr_I(R) = 0$. It can also be seen that $\M I^{n+1} : I = \M I^n$ for
all $n \geq 1$. Therefore $\depth F(I) > 0$.
\end{example}

\begin{example}\label{ex2}
Let $A = K[\![t^6, t^{11}, t^{15}, t^{31}]\!]$, $I = (t^6, t^{11},
t^{31})$ and $J = (t^6)$. Then, it can easily be verified that
$\ell(I^2/JI) = 1$ and $I^3 = JI^2$. Since $I^2 \cap J = JI$,
$G(I)$ is Cohen-Macaulay. It can also be seen that $t^{37} \in \M
I^2$, but $t^{37} \notin \M JI$. Therefore $F(I)$ is
not Cohen-Macaulay.
\end{example}

\end{document}